\documentclass[12pt]{amsart}
\usepackage[centertags]{amsmath} 
\usepackage{amsfonts} 
\usepackage{amscd} 
\usepackage{latexsym} 
\usepackage{amssymb} 
\usepackage{mathrsfs} 
\usepackage{euscript}

\newcommand{\B}[1]{{\mathbb #1}} 
 
\newcommand{\C}[1]{{\EuScript #1}}

\newcommand\R{\B R}

\newtheorem{theorem}[subsection]{Theorem}

\newtheorem{thm}[subsection]{Theorem} 
 
\newtheorem{corollary}[subsection]{Corollary} 
 
\newtheorem{cory}[subsection]{Corollary} 
 
\newtheorem{lemma}[subsection]{Lemma} 
 
\newtheorem{proposition}[subsection]{Proposition} 
 
\newtheorem{prop}[subsection]{Proposition} 
 
\newtheorem{question}[subsection]{Question}

\numberwithin{equation}{section}

\theoremstyle{definition} 
 
\newtheorem{definition}[subsection]{Definition}

\theoremstyle{remark} 
 
\newtheorem{remark}[subsection]{Remark} 
 
\newtheorem{remarks}[subsection]{Remarks}

\newcommand\theoref{Theorem~\ref} 
 
\newcommand\lemref{Lemma~\ref} 
 
\newcommand\propref{Proposition~\ref}

\numberwithin{figure}{section} 
 
\numberwithin{table}{section} 
 
\newcommand{\p}{{\partial}} 
 
\newcommand{\al}{{\alpha}}

\newcommand{\Om}{{\Omega}} 
 
\newcommand{\om}{{\omega}}

\newcommand{\ga}{{\gamma}} 
 
\newcommand{\Ga}{{\Gamma}}

\newcommand{\la}{{\lambda}}

\newcommand{\Si}{{\Sigma}}

\newcommand{\gf}{{\varphi}}

\newcommand{\Mo}{(M,\omega )}

\newcommand\Symp{\operatorname{Symp}} 
 
\newcommand\Ham{\operatorname{Ham}}

\def\Im{\operatorname{Im}}

\newcommand\id{\operatorname{Id}} 
 
\newcommand\cat{\operatorname{cat}} 
 \newcommand\qcat{\operatorname{Qcat}} 
\newcommand\wgt{\operatorname{wgt}} 
 
\newcommand\grad{\operatorname{grad}} 
 
\newcommand\sgrad{\operatorname{sgrad}} 
 
\newcommand\crit{\operatorname{crit}} 
 
\newcommand\Crit{\operatorname{Crit}}

\newcommand{\ms}{{\medskip}}

\def\la{\langle} 
 
\def\ra{\rangle} 
 
\def\ov{\overline} 
 
\begin{document}

\title{Symplectically aspherical manifolds} 
 
\author{Jarek K\c edra} 
 
\address{Mathematical Sciences, University of Aberdeen,  
  Meston Building, Aberdeen, AB243UE, Scotland, UK, and  
       Institute of Mathematics, University of Szczecin, 
Wielkopolska 15, 70-451 Szczecin, Poland} 
 
\email{kedra@maths.abdn.ac.uk} 
 
\urladdr{http://www.maths.abdn.ac.uk/\~{}kedra}

\author{Yuli Rudyak} 
 
\address{Department of Mathematics, University of Florida, 
358 Little Hall, Gainesville, FL 32601, USA} 
 
\email{rudyak@math.ufl.edu} 
 
\urladdr{http://www.math.ufl.edu/\~{}rudyak}

\author{Aleksy Tralle} 
 
\address{Department of Mathematics, University of Warmia and 
Mazury, 10-561 Olsztyn, Poland and  
Mathematical Institute, Polish Academy of Science, \'Sniadeckich 8, 00-950 
Warsaw,  
Poland} 
 
\email{tralle@matman.uwm.edu.pl} 
 
\urladdr{http://wmii.uwm.edu.pl/\~{}tralle}

 
 
 
\subjclass[2000]{Primary 53; Secondary 57}

\begin{abstract} 
 
This is a survey article on symplectically aspherical manifolds.

\end{abstract}

\maketitle

\section{Introduction}\label{S:intro} 
 
A symplectic form $\omega$ on a smooth manifold $M$ is called {\it  
symplectically  
aspherical}, if for any smooth map $f: S^2\to M$ one has the equality 
\begin{equation}\label{eq:sas} 
\int_{S^2}f^*\omega=0. 
\end{equation} 
The latter condition is often written in the form 
$ 
\omega|_{\pi_2(M)}=0. 
$ 
It is also useful to introduce a vector space  
$$ 
\Pi(M)=\Im\{h: \pi_2(M)\to H_2(M)\}\otimes \R, 
$$  
where $h$ is the Hurewicz map. In this notation, the symplectic asphericity 
condition can be written as 
$$ 
[\omega]|_{\Pi(M)}=0, \quad \text{ or }\la [\omega],x\ra=0 \text{ for all }  
x\in\Pi(M) 
$$ 
where $[\omega]$ is the (de Rham) cohomology class of the form $\omega$ and  
$\la-,-\ra$ is the Kronecker pairing. 
 
A symplectically aspherical manifold is, by definition, a manifold that is  
equipped with 
a symplectically aspherical form.  Such manifolds were originally introduced  
and used  
by Floer~\cite{F} in order to attack the Arnold conjecture. The usefulness of  
symplectically aspherical manifolds comes from the theory of $J$-holomorphic  
curves  
and the Floer  
homology, which are easier for symplectically aspherical manifolds because of  
the  
absence of  
bubbling effects (see \cite{F, H, HZ, MS1, MS2}). Since condition \eqref{eq:sas} 
is  
often  
imposed in many classical formulations of theorems in symplectic  
topology~\cite{LO, O}, 
it is  
worthwhile to describe this class of manifolds. In the last decade a substantial  
understanding of  
symplectically aspherical manifolds was gained in \cite{G1, IKRT, KRT, R2, RO,  
RT1, S}.  
However, there are still open interesting questions, and the whole subject  
becomes a  
rich mathematical theory involving many important topological and symplectic  
techniques.  
Motivated by this, we present a survey on recent developments in the theory of 
symplectically aspherical manifolds. We emphasize, however, that we are 
mainly concerned with the ``soft'' symplectic geometry, where the tools of 
algebraic 
topology are applicable. We use hard analytical tools ( the pseudoholomorphic 
curves  in 
Section \ref{S:arnold}, the applications of the action spectrum in Section 
\ref{S:applications}) rather than develop them. 
 
The main topics described in this article are the 
following. 
\begin{enumerate} 
\item Elementary topological properties of symplectically aspherical manifolds; 
\item Lusternik--Schnirelmann category of closed symplectically aspherical 
manifolds; 
\item Applications of symplectic asphericity to the classical Arnold conjecture  
on  
fixed points of symplectic diffeomorphisms; 
\item Constructions of symplectically aspherical manifolds; 
\item Discussion on fundamental groups of closed symplectically aspherical  
manifolds; 
\item Applications of symplectic asphericity to the theory of the group of  
symplectic  
diffeomorphisms and to the symplectic action spectrum; 
\item Discussion on symplectically hyperbolic manifolds; 
\item Applications of symplectic asphericity to Lie group actions on closed  
symplectic  
manifolds. 
\end{enumerate} 
 
We survey the known results and discuss research perspectives and conjectures.   
     
Throughout this paper, all manifolds are assumed to be closed and connected  
(all exceptions are explicitly stated).

\section{Preliminaries on symplectically aspherical manifolds}\label{S: prelim} 
 
In this section we present some facts which we use in the sequel.   
Throughout  
the section the manifolds are not assumed to be closed.

\begin{prop} 
\label{induced} 
 
Let $\omega$ be a symplectically aspherical form on a manifold $N$, and let  
$g: M \to N$ be a map such that $g^*\omega$ is a symplectic form on a manifold  
$M$. Then $(M, g^*\omega)$ is symplectically aspherical. In particular, a 
covering manifold  
over a symplectically aspherical manifold is symplectically aspherical. 
\end{prop} 
 
\begin{proof} Since the form $\omega$ on $N$ is symplectically aspherical, we  
conclude  
that 
$$ 
\langle g^*\omega, h(a)\rangle=\langle\omega, g_*h(a)\rangle= 
\langle\omega, h(g_*(a))\rangle=0 
$$ 
for all $a\in \pi_2(M)$. Thus, $g^*\omega$ is symplectically aspherical. 
\end{proof}

\ms We need the following homotopic characterization of symplectically 
aspherical manifolds. Given a group $\pi$, recall that the Eilenberg--Mac Lane  
space $K(\pi,1)$ is  a connected $CW$-complex with fundamental  
group $\pi$ and such that $\pi_i(K(\pi,1))=0$ for $i>1$. It is well known  that  
the homotopy type of such space is completely determined by
$\pi$. Moreover, for every connected $CW$-space $X$ with $\pi_1(X)=\pi$ there
exists a map $f: X\to K(\pi,1)$ that induces an isomorphism on the
fundamental groups, and this map is unique up to homotopy. To
construct such a map $f$, attach to $X$ cells of dimensions $>2$ in
order to kill all the higher homotopy groups. The resulting space is
$K(\pi,1)$ and $f$ is the inclusion of $X$.

\begin{prop}[{\rm  \cite[Lemma 4.2]{LO}, \cite[Corollary 2.2]{RT}}] 
\label{p:aspher} 
 
Let $(M,\omega)$ be a symplectic manifold.  
The following three  
conditions are equivalent: 
 
\par {\rm (i)} the form $\omega$ is symplectically aspherical;  
 
\par {\rm (ii)} if a map $f: M\to K(\pi_1(M),1)$ induces an  
isomorphism on the fundamental groups, then 
$$  
[\omega]\in\Im\,\{f^*: H^2(K;\R)\to H^2(M;\R)\}; 
$$ 
 
\par {\rm (iii)} there exist a group $\tau$ and a map $g: M\to
K(\tau,1)$ such that  
$$  
[\omega]\in\Im\,\{f^*: H^2(K(\tau,1);\R)\to H^2(M;\R)\}. 
$$ 
\end{prop}

\ms In the context of symplectic asphericity, it seems natural to ask if there  
is a manifold $M$ that possesses two symplectic forms $\omega_1$ and $\omega_2$  
such that $\omega_1|_{\pi_2(M)}=0$ and $\omega_2|_{\pi_2(M)} \ne 0$. The 
following  
Proposition answers this question. 
 
\begin{prop} Suppose that a closed manifold $M$ admits a symplectic form. Then 
the  
following two conditions are equivalent: 
 
\par{\rm (i)} $\Pi(M)=0$; 
 
\par{\rm (ii)} Every symplectic form on $M$ is symplectically aspherical. 
\end{prop}

\begin{proof} Only the implication (ii) $\Longrightarrow $(i) needs a proof.  
So, assume  
that $\Pi(M)\ne 0$ and consider a symplectic form $\omega$ on $M$. If 
$\omega|_{\Pi}\ne 0$ then we are done. So, assume that $\omega|_{\Pi}= 0$. Since 
$\Pi(M)\ne 0$,  
there exists a closed 2-form $\sigma$ with $\la[\sigma], x\ra \ne 0$ for some 
$x\in  
\Pi(M)$, i.e. $\sigma|_{\pi_2(M)}\ne 0$. Now, the form $\ga=\omega +  
\lambda\sigma$  
is symplectic for $\lambda$ small enough, and $\ga|_{\pi_2(M)}\ne 0$. 
\end{proof}  
 
\section{Lusternik--Schnirelmann category of symplectically aspherical  
manifolds} 
\label{S:LS}

In this section we prove that, for any (closed) symplectically aspherical  
manifold $M^{2n}$, the number of critical points of any smooth function $f: M 
\to \R$ is  
at least $2n+1$. We also describe some results that will be used in our 
discussion of  
the Arnold conjecture in Section 4.

\subsection{Lusternik--Schnirelmann Theorem} 
 
Given a smooth function $F: M \to \R$ on a smooth manifold $M$, let $\crit F$  
denote  
the number of critical points of $F$. Put $\Crit M=\min\{\crit F\}$, where the 
minimum   
runs  
over all smooth functions  $F: M \to \R$

\begin{definition}[\cite{Fox}] 
 
\label{def:cat-map} 
 
Let~$f: X \to Y$ be a map. The 
{\it Lusternik--Schnirelmann category of~$f$}, denoted~$\cat(f)$, is 
defined to be the minimal integer~$k$ such that there exists an open 
covering~$\{U_0, \ldots, U_k\}$ of~$X$ with the property that each of 
the restrictions~$f|U_i\colon U_i \to Y$,~$i=0,1, \ldots, k$ is 
null-homotopic. If such a covering does not exist we say that $\cat(f)$ is not  
defined. 
 
The {\it Lusternik--Schnirelmann category~$\cat X$ of a space~$X$} is 
defined as the category~$\cat (\id_X)$ of the identity map. 
\end{definition}

For the proof of the following Lusternik--Schnirelmann Theorem, see~\cite{CLOT}. 
 
\begin{theorem}\label{t:ls} 
For every closed manifold $M$ we have the inequality $\Crit M \ge \cat M+1$. 
\end{theorem} 
 
This theorem admits the following generalization. Given a flow  
$\{\gf_t: X \to X, t\in \R\}$,  
a rest point of the flow is defined as a point $x\in X$ such that $\gf_t(x)=x$  
for  
all $t\in \R$. The flow is {\it gradient-like} if there exists a function $F: X  
\to R$  
(called a {\it Lyapunov function}) such that  
$F(\gf_t(x))>F(\gf_s(x))$ whenever $t<s$ and $x$ is not a rest point of the 
flow. 
 
\begin{thm}\label{t:flow} 
Let $X$ be a compact space and $f: X \to Y$ be a map such that $\cat(f)$ is  
defined. Let $\{\gf_t\}$ be a gradient-like flow on a compact space $X$. Then  
the number of rest points of the flow is at least $1+\cat(f)$. 
\end{thm} 
 
The proof can be found in~\cite{CLOT}. Note that \theoref{t:flow} implies  
\theoref{t:ls}. Indeed, given a closed manifold $M$, let $f$  be the identity  
map in \theoref{t:flow}. Now, given a smooth function $F: M \to R$, we have the  
gradient-like flow $-\grad F$ (this explains the name gradient-like),  
whose rest points are exactly the critical points of $F$.    
 
\subsection{Category weight} 
The following definition is a homotopy invariant version of a construction of  
Fadell--Husseini~\cite{FH}. It was suggested by Rudyak~\cite{R1} and 
Strom~\cite{St}.

\begin{definition}\label{def:swgt} 
The \emph{category weight}~$\wgt(u)$ of a non-zero cohomology class~$u \in 
H^*(X; A)$ is defined as follows: 
\begin{equation*} 
\label{31} 
\wgt(u)\ge k \Longleftrightarrow \{\gf^*(u)=0 {\rm\ for\ every\ } \gf\colon F 
\to X 
{\rm\ with\ } \cat(\gf) < k\}. 
\end{equation*} 
\end{definition} 
 
\begin{prop}[\cite{R1,St,CLOT}] 
\label{p:wgtprops} 
Let $A$ denote a coefficient ring. Category weight has the following properties. 
\begin{enumerate} 
\item~$1\le \wgt(u) \leq \cat(X)$, for all~$u \in \widetilde 
H^*(X;A), u\ne 0$. 
\vskip3pt 
\item For every~$f\colon Y \to X$ and~$u\in H^*(X;A)$ with 
$f^*(u)\not = 0$ we have 
$\cat(f) \geq \wgt(u)$ and~$\wgt(f^*(u)) \geq \wgt(u)$. 
\vskip3pt 
\item For~$u\in H^*(X;A)$ and~$v\in H^*(X;A)$ we have 
\[ 
\wgt(u\cup v) \geq \wgt(u) + \wgt(v). 
\]  \vskip3pt 
\item For every~$u \in H^s(K(\pi,1);A)$,~$u\ne 0$, we have 
$\wgt(u)=s$.  \vskip3pt 
\item For every~$u \in H^s(X;A)$,~$u\ne 0$, we have 
$\wgt(u)\le s$.  \vskip3pt 
\end{enumerate} 
\end{prop} 
 
\subsection{Symplectic asphericity input} 
Here we show the effect of symplectic asphericity on the category  
weight of the symplectic class $[\omega]$. 
 
\begin{thm}[\cite{RO}] 
\label{t:omega} 
If $(M,\omega)$ is a symplectically aspherical manifold then $\wgt([\omega])=2$. 
\end{thm} 
 
\begin{proof} 
It follows from \propref{p:aspher} that $[\omega]=f^*a$ for some $a\in  
H^2(K(\pi_1(M),1);\R)$  
and some $f: M \to K(\pi_1(M),1)$. Now, $\wgt(a)=2$ by item (4) of  
\propref{p:wgtprops}, and thus  
$\wgt([\omega])=2$ by items (2) and (5) of the same proposition.  
\end{proof} 
 
\begin{cory}\label{c:lscat} 
If $(M^{2n},\omega)$ is a symplectically aspherical manifold then $\cat M=2n$ 
and  
$\Crit M =2n+1$.  
\end{cory} 
 
\begin{proof}  Since $[\omega]^n\ne 0$, we conclude that  
$\wgt([\omega]^n)=2n$ by \theoref{t:omega} and items (3) and (5) of  
\propref{p:wgtprops}. Furthermore, $\cat M \ge \wgt([\omega]^n)=2n$, and so  
$\Crit M \ge 2n+1$.  Finally, according to Takens~\cite{T}, we have $\Crit N \le 
\dim N+1$ for every closed  
connected manifold $N$, and the result follows. 
\end{proof} 
 
\section{The Arnold conjecture}\label{S:arnold} 
 
Let $(M, \omega)$ be a symplectic manifold. It is well-known that there exist a  
Riemannian metric $g$ and and almost complex structure $J$ on $M$ such that  
$\om (\xi, J\eta)=g(\xi, \eta)$.  
 
\subsection{Hamiltonian diffeomorphisms}  
Given a function $F: M \to R$, define a symplectic-gradient vector field  
$\sgrad F$ (frequently denoted also by $X_F$) by the condition 
\begin{equation} 
\omega(\sgrad F, \xi)=-dF(\xi) 
\end{equation} 
for all vector fields $\xi$. It is easy to see that $\sgrad F=J\grad F$ where  
$\grad F$ is  
taken with respect to the metric $g$. 
 
Now, consider a smooth function $H: S^1\times M \to \R$ and put  
$H_t(x)=H(t,x)$.  
Consider the non-autonomic differential equation 
\begin{equation}\label{eq:ham} 
\dot x(t)=\sgrad H_t(x(t))  
\end{equation} 
This equation yields a time-dependent flow $\Psi=\{\phi_t\}=\{\phi_t^H\}$ on  
$M$. Namely, if $x(t)$ is a solution of \eqref{eq:ham} with $x(0)=p\in M$, then  
$\phi_t(p)=x(t), t\in \R$.  
 
\begin{definition} 
A diffeomorphism $\phi: M \to M$ is {\it Hamiltonian} if there exists a  
function (which is called a Hamiltonian) 
$H:S^1\times M \to \R$  such that  
$\phi=\phi_1^H$. We also say that $\phi$ is a time-1 map of the Hamiltonian $H$. 
The set of all Hamiltonian diffeomorphism is denoted by $\Ham\Mo$. 
\end{definition} 
 
\noindent {\bf The Arnold conjecture.}  {\it For every Hamiltonian 
diffeomorphism  
$\phi: M\to M$, the number of its fixed points is at least $\Crit M$.} 
 
\begin{remark} 
There are several versions of the Arnold conjecture. The one above 
is the closest to the ``soft'' side of  symplectic geometry. 
If we assume that a Hamiltonian diffeomorphism is non-degenerate, 
that is its graph intersects the diagonal transversely,  
then the conjecture claims that 
the number of fixed points is at least the number of critical 
points of any Morse function. This conjecture is neither proved nor disproved 
yet, even in the symplectically aspherical case. We mention here the 
papers~\cite{FO,LT} where the number of fixed points is estimated from below by 
the sum of Betti numbers.
\end{remark}

\subsection{Floer's approach to the Arnold conjecture} Floer~\cite{F}      
suggested the following way to attack  the Arnold conjecture for  
symplectically aspherical manifolds. If $p$ is a fixed point of a Hamiltonian  
symplectomorphism $\phi=\phi_1^H$ then $\phi_t(p), t\in \R$ is a 1-periodic  
orbit, i.e. a loop $x: S^1 \to M$. Hence, we can count fixed points of $\phi$  
by counting 1-periodic solutions of equation \eqref{eq:ham}. So, we can try to  
pose a variational problem on the loop space of $M$ whose solutions (extremals,  
critical loops) are the 1-periodic solutions of \eqref{eq:ham}. Then we can  
apply the Lusternik--Schnirelmann theory to estimate of the number of  
extremals. Note that we will count only contractible 1-periodic orbits. Let us 
explain this in a bit more detail. 
 
\subsection{Variational reduction} Let $H: S^1\times M \to \R$ be a Hamiltonian  
on a symplectically aspherical manifold $\Mo$. Given a contractible smooth loop  
$x: S^1\to M$, we set 
\begin{equation}\label{eq:action}  
\C A_H(x)=\int_{D^2}y^*\om-\int_0^1H(t, x(t))dt 
\end{equation} 
where $y: D^2\to M$ is an extension of $x$. We call the functional $\C A_H$ on  
contractible loops  the {\it symplectic action}. Note that the symplectic  
action is well-defined (it does not depend on the extension $y$) because of   
symplectic asphericity. So, we have the map 
$$ 
\C A_H: C^{\infty}_c(S^1,M):=\{\text{contractible smooth maps } S^1 \to M\}  
\longrightarrow \R. 
$$     
If we take the derivative of  $\C A_H$ in the direction of the vector field 
$\xi$ along  
$x(S^1)$ 
(regarded as a tangent vector to a loop $x$), we get 
\begin{equation} 
\begin{aligned} 
D \C A_H(x)(\xi)&=\int_{S^1}\om_{x(t)}(\dot x(t), \xi)- dH_t(x(t)) \xi(t)dt\\ 
&=\int_{S^1}\om_{x(t)}(\dot x(t)-\sgrad H_t(x(t)), \xi))dt. 
\end{aligned} 
\end{equation} 
So, if $x=x(t)$ is a critical orbit of $\C A_H$, that is $D \C A_H(x)(\xi)=0$ 
for  
all $\xi$,  then $\dot x(t)-\sgrad H_t(x(t))=0$, (i.e. $x(t)$ is a 1-periodic  
solution of equation~\eqref{eq:ham}).

In order to proceed, we must consider the ``gradient flow'' of $\C A_H$.  
However, here we have many analytical difficulties that do not allow us to  
construct the gradient flows directly, cf. ~\cite[Section 6.5]{HZ}.  Floer  
regards the gradient flow lines as maps  
\begin{equation}\label{eq:u} 
u: \R \times S^1\to M,\quad (s,t)\mapsto u(s,t)=u(s,t+1) 
\end{equation} 
such that 
\begin{equation}\label{eq:flowline} 
\frac{\p u}{\p s}+J(u)\frac{\p u}{\p t}+\grad H_t(u)=0. 
\end{equation} 
Floer~\cite{F} then obtained the following variational reduction theorem. 
 
\begin{thm}\label{t:reduction} 
Let $(M,\om)$ be a symplectically aspherical manifold. Assume also that  
$c_1(M)$ vanishes on $\pi_2(M)$. Let $\phi: M \to M$ be a Hamiltonian  
symplectomorphism. Then there exist a map $f: X \to M$ and a gradient-like  
flow $\Phi$ on $X$ with the following properties:  
\begin{enumerate} 
\item~$X$ is a compact metric space. 
\vskip3pt 
\item~The number of fixed points of $\phi$ is bounded from below by the number  
of rest points of $\Phi$. 
\vskip3pt 
\item The map $f^*: H^*(M;A) \to H^*(X;A)$ is a monomorphism for any  
coefficient group $A$. 
\end{enumerate} 
\end{thm}

Now we describe the data $X,f$ and $\Phi$ of \theoref{t:reduction}. Given a map  
$u: \R \times S^1\to M$ and $\tau\in \R$, we define $u(\tau): S^1\to \R$, by  
$u(\tau)(s)=u(\tau,s)$.  
 
Let $X$ be a space of smooth maps $u: \R \times S^1$ such that 
\begin{enumerate} 
\item~$u$ satisfies the equation \eqref{eq:flowline}. 
\vskip3pt 
\item $u(0)$ (and hence $u(\tau)$ for all $\tau$) is a contractible loop. 
\vskip3pt 
\item The function $\gf: \R \to \R,\, \gf(s)=\C A_H(u(s))$ is bounded. 
\end{enumerate} 
It turns out  that $X$ is compact. It is worth  mentioning that the proof  
of compactness uses the symplectic asphericity of $\Mo$. 
 
\ms 
We define a flow $\Phi=\{\gf_{\tau}, \tau\in \R\}$ on $X$ by setting  
$\gf_{\tau}(u(s,t))=u(s+\tau, t)$. It can be proved that the flow is  
gradient-like with associated Lyapunov function  
$$ 
F: X \to \R,\quad F(u)=\C A_H(u(0)). 
$$ 
Furthermore, if $u\in X$ is a rest point of $\Phi$ then $\p u /\p s=0$. So, if  
we put $x=u(0): S^1\to M$ then 
$$ 
J(x)\frac{dx}{dt}+\grad H_t(x)=0, 
$$ 
or 
$$ 
\frac{dx}{dt}=J\grad H_t(x)=\sgrad H_t(x). 
$$ 
Note that the latter equation is obtained by applying $J$ to both sides of the 
preceding   
equation and using $J^2=-I$. 
So, $x$ is a 1-periodic solution of the equation \eqref{eq:ham}, and therefore  
$x(0)$  
is a fixed point of the Hamiltonian diffeomorphism $\phi$. Thus, the number of  
fixed  
points of $\phi$ is at least the number of rest points of $\Phi$. 
\ms 
Finally, we define a map $f: X \to M$ by setting $f(u)=u(0,0)$. The proof of  
the monomorphicity of $f^*$ is difficult. 
\begin{remarks} 
1. The above description of the space $X$ is taken from the book \cite{HZ}.  
In his original paper~\cite{F} Floer obtained the space $X$ as a certain space 
of  
contractible loops $S^1\to M$. As we have seen, critical points of $\C A_H$ 
are 1-periodic solutions of the equation \eqref{eq:ham}. 
Generally, if one has a flow on a space, say,  
$Y$, one can consider a new space $\ov Y$ whose points are the flow lines on  
$Y$,  
and define the flow on $\ov Y$ via the time-shift. This describes the passages  
from Floer's interpretation to that of Hofer and Zehnder.  
 
2. It seems that the condition $c_1(M)|_{\Pi(M)}=0$ in \theoref{t:reduction} is  
redundant, cf.~\cite[Remark on p. 250]{HZ}, but as far as we know, nobody has 
written  
this down yet in the literature. 
\end{remarks} 
 
\subsection{Proof of the Arnold conjecture for symplectically aspherical  
manifolds} Now, basing  our argument on  \theoref{t:reduction} and results of  
Section~\ref{S:LS}, we prove the Arnold conjecture under assumptions of  
\theoref{t:reduction}. 
 
\begin{thm}[\cite{R2, RO}] 
The Arnold conjecture holds for symplectically aspherical manifolds with  
$c_1(M)|_{\Pi(M)}=0$. 
\end{thm} 
 
\begin{proof} 
We use the notation from \theoref{t:reduction}. It suffices to prove that the  
number of rest points of $\Phi$ is at least $2n+1=\Crit M$. Since  
$f^*([\omega]^n)\ne 0$, we conclude that $\cat f \ge \wgt [\om]^n=2n$ by  
\propref{p:wgtprops} and \theoref{t:omega}. So, by \theoref{t:flow}, the number  
of rest points of $\Phi$ is at least  
$2n+1=\Crit M$. 
\end{proof} 
 
 \subsection{Lagrangian submanifolds and their intersections} A Lagrangian 
submanifold of a symplectic manifold $(V^{2n},\omega)$ is a smooth submanifold 
$L^n$ of $V$ such that $\omega|_{L}=0$. Given a Hamiltonian diffeomorphism 
$\psi: 
V \to V$, the question on the number $\#(\psi(L)\cap L)$ of intersection points 
of $\psi(L)$ and $L$ can be considered as a generalization of the Arnold 
conjecture. 
Indeed, given a symplectic manifold $\Mo$, the diagonal $M$ of the symplectic 
manifold 
$(M\times M, \omega \times (-\omega))$ is Lagrangian. Furthermore, given a 
Hamiltonian 
diffeomorphism $\phi: M \to M$, define $\psi: M \times M \to M \times M$ as 
$\psi(x,y)=(\phi(x),y)$. Then the number $\#(\psi(M)\cap M)$ is exactly the 
number of fixed points of $\phi$. 

There is a large literature on Lagrangian intersections, we 
mention~\cite{EG,F,H}. In 
this context the symplectic asphericity appeared in~\cite{F,H} in the form 
$\pi_2(V,L)=0$.

It is well known that the total space $T^*L$ of the cotangent bundle of a smooth 
manifold 
$L$ possesses a canonical symplectic form $\omega$. Furthermore, the zero 
section $L$ is a 
Lagrangian submanifold of $(T^*L,\omega)$. Moyaux and Vandembroucq~\cite{MV} 
estimated from 
below the number $\#(\psi(L)\cap L)$ where $\psi: T^*L \to T^*L$ is a 
Hamiltonian 
diffeomorphism with a compact support. Generally this number is bounded by a 
certain numerical invariant $\qcat L$, but if $L$ is symplectically aspherical 
then 
$\#(\psi(L)\cap L)\ge \Crit L$,~\cite[8.5.1]{CLOT}.

\section{Constructions of symplectically aspherical manifolds} 
\label{S:constructions} 
 
In this section we present various constructions of symplectically 
aspherical manifolds. Mostly, they are based on the known constructions 
of symplectic manifolds  which yield the symplectically aspherical property 
under additional hypotheses.

\subsection{Branched coverings}\label{SS:gompf} 
Here we follow \cite{G2}. Take a symplectically aspherical manifold $(X,\om_X)$,  
choose a symplectic submanifold 
$B\subset X$ of codimension 2 and construct $\Mo$ 
as a covering of $X$ branched along $B$. By construction 
the class $[\om]$ of the symplectic form is the pull-back 
of the class $[\om_X]$, and so $\omega$ is symplectically aspherical by 
\propref{induced}.   
Note that Proposition 2.1 applies since outside the branch locus, we have a true 
covering. 
 
By choosing the branching locus in a clever way, Gompf  
constructed the first examples of symplectically aspherical 
manifolds with non-trivial second homotopy group. Moreover, 
he proved that the symplectic asphericity and the vanishing 
of the first Chern class on spheres are independent conditions~\cite[Theorem 
7]{G2}).  
Also, Gompf constructed both K\"ahler and non-K\"ahler 
examples. 
 
\subsection{Homogeneous spaces}\label{SS:homogeneous} 
 
Let $G$ be a simply connected solvable group and $\Ga \subset G$ 
a uniform lattice. Recall that a uniform lattice is a discrete 
subgroup such that the quotient $G/\Ga$ is compact. Since 
a simply connected solvable group is contractible (in fact 
diffeomorphic to $\B R^n$) the quotient is a $K(\Ga,1)$-space. 
If $G$ admits a $\Ga$-invariant symplectic form then the 
quotient $G/\Ga$ is a symplectic manifold. This happens in 
many cases. 
 
Let $\text{Lie}(G)$ denote the Lie algebra of a Lie group $G$. 
A group $G$ is called completely solvable if the eigenvalues 
of the operators $\text{ad}_X:\text{Lie}(G)\to \text{Lie}(G)$ are 
real for any $X\in \text{Lie}(G)$. In this case the Hattori 
theorem \cite{Ha} states that 
$$ 
H^*(G/\Gamma;\B R)\cong H^*(\Lambda^*(\text{Lie}(G));\B R), 
$$ 
where the right-hand side is the cohomology of the Lie algebra 
of $G$. In particular, every cohomology class in $H^*(G/\Ga; \R)$ can be 
represented by a left-invariant form.  
Hence in order to show that $G/\Gamma$ is symplectic 
it is enough to find a cohomology class $a\in H^2(G/\Gamma;\B R)$ 
such that its top power is non-zero. A manifold admitting 
such class $a$ is called cohomologically symplectic. We summarize 
the above discussion in the following. 
 
\begin{theorem}[{\cite[Lemma 4.2]{IKRT}}] 
\label{T:homogeneous} 
Let $G$ be a simply connected completely solvable Lie group and 
$\Gamma \subset G$ a  uniform lattice. 
A homogeneous space $M=G/\Gamma$ is symplectically aspherical 
if and only if it is cohomologically symplectic. \qed 
\end{theorem} 
 
\subsection{Symplectic bundles}\label{SS:fibrations} 
 
A locally trivial bundle $\Mo \to E\to B$ is called 
symplectic if its structure group is the group of 
symplectic diffeomorphisms of $\Mo$. The following theorem 
gives conditions implying that the total space of a 
symplectic bundle is symplectically aspherical. 
 
\begin{theorem}[{\cite[Theorem 7.4]{IKRT}}] 
\label{T:fibrations} 
Let $\Mo$ and $(B,\omega_B)$ be symplectically 
aspherical manifolds. The total space of a symplectic 
bundle $\Mo\stackrel{i}\to E\stackrel{p}\to B$ is 
symplectically aspherical if the following conditions 
are satisfied: 
\begin{enumerate} 
\item 
there exists $\Om\in H^2(E)$ such that $i^*\Om=[\om]$ and 
\item 
the class $[\om]$ vanishes on $(i_*)^{-1}(\Pi(E))$. 
\end{enumerate} 
\end{theorem} 
 
The first condition ensures the existence of a symplectic 
form on the total space. This was proved by Thurston in 1976, see~\cite[Theorem 
6.3]{MS1}. 
 The second condition ensures the symplectic asphericity of the symplectic form  
 coming from the Thurston   
construction.  
 
\subsection{Symplectic surgery}\label{SS:surgery} 
 
In \cite{G1}, Gompf proved that certain surgery can 
be performed symplectically. More precisely, let 
$j_i:(N,\omega_N)\to (M_i,\omega_i)$, $i=1,2$, be disjoint symplectic 
embeddings of codimension two. Suppose that their normal bundles  
have opposite Euler classes. Cut out small tubular 
neighborhood of the images $j_i(N)$ and glue the  
remaining part along the cutting locus to obtain a new 
manifold $X=M_1\cup_N M_2$. This manifold admits a symplectic form 
which is equal to $\om$ away from the gluing locus. 
The next theorem gives a condition under which this 
construction produces an aspherical symplectic form. 
 
\begin{theorem}[{\cite[Theorem 6.3]{IKRT}}] 
\label{T:surgery} 
Let $(M_1,\om_1),(M_2,\om_2),(N,\om_N)$ be  
symplectic manifolds and $j_i:N\to M_i$ symplectic  
embeddings with opposite normal bundles.  
If $M_1$ and $M_2$ are symplectically aspherical 
and $(j_i)_*:\pi_1(N)\to \pi_1(M_i)$ are monomorphisms, 
then the Gompf symplectic sum 
$M_1\cup_N M_2$ is symplectically 
aspherical. 
\end{theorem} 
 
\subsection{Lefschetz fibrations}\label{SS:lefschetz} 
An excellent exposition on 
Lefschetz pencils and fibrations can be found in 
\cite{GS}.   
Let $X$ be a compact, connected, oriented, smooth 4-manifold, possibly with 
boundary. {\it A Lefschetz fibration structure} on $X$ is a surjective map $f: 
X\rightarrow \Sigma$, where $\Sigma$ is a compact, connected, oriented surface 
and $f^{-1}(\partial\Sigma)=\partial X$. Furthermore, the following is required:  
\begin{itemize}  
\item the set $\{q_1,...,q_n\}$ of critical points of $f$ is finite;  
\item $f(q_i)\not=f(q_j)$ for $i\not=j$;  
\item if $b\in\Sigma$ is a regular value of $f$ then $f^{-1}(b)$ is a closed 
connected orientable surface;  
\item there exist an orientation-preserving complex charts $\varphi_i: 
U_i\rightarrow \mathbb{C}^2$ with   
$q_i\in U_i\subset X$ and $\psi_i: V_i\rightarrow \mathbb{C}$ with $f(q_i)\in 
V_i\subset f(U_i)\subset \Sigma$ such that   
$\psi_i\circ f\circ\varphi^{-1}_i:\varphi(U_i)\rightarrow\mathbb{C}$ has the 
form $(x,y)\rightarrow x^2+y^2$.  
\end{itemize}   
It is a celebrated result of Donaldson \cite{D1} that 
a closed symplectic 4-manifold admits a structure of 
an oriented Lefschetz pencil, and, therefore, becomes a Lefschetz fibration 
after blowing up in a finite number of points. On the other hand,  
according to Gompf and Thurston, an oriented Lefschetz  
fibration admits a symplectic structure. 
 
In principle, this gives a classification of 4-dimensional 
symplectic manifolds in terms of the monodromy of 
a Lefschetz pencil. One of the main results of \cite{KRT} 
provides a condition under which a Lefschetz fibration 
admits an aspherical symplectic form. In the sequel we denote by $\Sigma_g$ a 
closed orientable surface of genus $g$. 
 
\begin{theorem}[{\cite[Proposition 3.3]{KRT}}] 
\label{T:lefschetz} 
Let $F$ be a closed connected orientable surface.  
Let $F\to X \to \Si_g$ be a symplectic Lefschetz fibration 
such that the inclusion of the fiber induces a 
non-trivial map  
$$ 
H^2(X;\B R)\to H^2(F;\B R)=\B R. 
$$ 
Put $Y =F\times \Si_h$ with the product symplectic structure.  
If $g+h>0$ then the Gompf symplectic 
fiber sum $X \#_F Y$ is symplectically aspherical. 
\end{theorem} 
 
\section{Fundamental groups of symplectically aspherical 
manifolds}\label{S:group} 
 
\ms For brevity, we call a group $\Gamma$ {\it symplectically aspherical} if it  
can be realized as the fundamental group of a closed symplectically aspherical  
manifold. 
\begin{question}\label{Q:group} 
What groups are symplectically aspherical? 
\end{question} 
 
This question has various motivations. 
The first  one belongs to a class of questions revolving around whether  a given 
geometric structure  
imposes  
restrictions on the algebraic topology of the underlying 
manifold. The second motivation comes from the problem of describing  properties 
of the fundamental group which  
determine the geometry of the manifold. In this section,  
we shall present  results in the first direction. 
 
An example of the second approach is Corollary \ref{C:hyperbolic}. 
 
\ms 
Questions similar to 6.1 are still unanswered in the case of complex 
projective or K\"ahler manifolds \cite{ABCKT}. It easily follows 
from the Lefschetz property that the first Betti number of the 
fundamental group of a closed K\"ahler manifold is even. 
According to Gromov and Shoen certain condition on the 
fundamental group of a closed K\"ahler manifold $M$ 
implies that it admits a  holomorphic 
mapping onto a Riemann surface. Another example is the 
the Shafarevich conjecture which  
states that if $M$ is a complex manifold with  $\pi_1(M)$ large, then the  
universal cover  
$\widetilde M$ is a Stein manifold. Here the fundamental group 
is called {\em large} if for any non-constant holomorphic  
map $f:X\to M$ the image $f_*(\pi_1(X))\subset \pi_1(M)$ 
is infinite. 
 
\ms 
Gompf proved in \cite{G1} that every finitely presented group can be realized as 
the fundamental group of a closed symplectic manifold. 
In contrast, we shall show that it is not the case for 
symplectically aspherical manifolds. For example, observe 
that symplectically aspherical groups have to be infinite. 
Indeed, since the symplectic asphericity is preserved by  
finite coverings, if $\Mo$ is symplectically 
aspherical and has finite fundamental group then 
its universal covering $(\widetilde M,\widetilde \om)$ 
is a symplectically aspherical simply connected closed manifold.  
This is impossible  because a non-zero class $[\widetilde \om]$  
 does not vanish on $H_2(\widetilde M;\B Z)=\pi_2(\widetilde M)$.  
 
\begin{proposition}\label{P:group} 
If a group $\Gamma$ is a fundamental group of a symplectically aspherical 
manifold then either 
 
\begin{enumerate} 
\item 
$\Gamma \cong \pi_1(\Si)$, where $\Si$ is a closed  
oriented surface, or 
\item 
there exists $\Om\in H^2(\Gamma;\B R)$ with $\Om^2\neq 0$. 
\end{enumerate} 
\end{proposition} 
 
\begin{proof} 
Let $\Mo$ be a symplectically aspherical. The case of 
dimension 2 is trivial so let us assume that $\dim M>2$. 
We know that $[\om]=c^*\Om$, where $c:M\to K(\pi_1(M),1)$ 
is the classifying map. Since $[\om]^2\neq 0$ we get that 
$\Om^2\neq 0$. 
\end{proof} 
 
In particular, no group of  real cohomological  
dimension 3 is symplectically aspherical.  
 
Constructions of symplectically aspherical groups are based on the 
constructions presented in Section \ref{S:constructions}. 
Using the Lefschetz fibrations it is possible to give 
a complete classification of symplectically aspherical 
Abelian groups. 
 
\begin{theorem}[{\cite[Theorem 1.2]{KRT}}] 
\label{T:abel} 
A finitely generated Abelian group $\Gamma $ is  
symplectically aspherical if and only if either $\Gamma \cong \mathbb Z^2$ or  
$\text{\em rank}(\Ga)\geq 4$. 
\end{theorem} 
 
\ms 
As we mentioned in the beginning of this section, 
the first Betti number of a K\"ahler 
group is even. The next result 
states that there is no such restrictions in the 
symplectically aspherical case. 
 
\begin{theorem} 
Any non-negative integer number that is different from one can be realized as 
the first Betti number of 
a symplectically aspherical manifold. 
\end{theorem} 
 
\begin{proof} 
{\bf The case $b_1=0$.} Let $\Ga$ be a uniform lattice in $SU(2,1)$ 
. It acts 
freely by isometries on the complex hyperbolic plane and hence 
preserves the K\"ahler structure on it. Thus the quotient is an 
aspherical closed K\"ahler (and hence symplectic) manifold. Finally, 
$b_1(\Gamma)=0$ \cite{G, Ma}. 
 
{\bf The case $b_1=2$.} The two dimensional torus serves as an example. 
 
{\bf The case $b_1=3$.} The Kodaira--Thurston manifold $KT$, see \cite{MS1} is 
an example. 
More precisely, $KT$ is the product $S^1 \times (N_3(\B R)/N_3(\B Z))$, 
where $N_3(\B K)$ denotes the upper triangular matrices with 
coefficients in $\B K$. Hence the fundamental group is 
isomorphic to $\B Z\oplus N_3(\B Z)$ and it is easy to 
see that its first Betti number is equal to 3. 
 
{\bf The case $b_1\ge 4$.} This follows from Theorem \ref{T:abel}, or you can 
consider the products of tori and Kodaira--Thurston manifolds. 
\end{proof} 
 
\begin{remark} 
Unfortunately, we are unable to construct a symplectically aspherical 
manifold with the first Betti number equal to one. Nevertheless, 
we believe that such manifolds exist. 
\end{remark}

\section{Further applications of symplectic asphericity} 
\label{S:applications} 
 
\subsection{The action spectrum} 
\label{SS:schwarz} 
 
Let $\Mo$ be a symplectic manifold and let $H:S^1\times M\to \B R$ be a  
Hamiltonian. 
A contractible solution $x:S^1\to M$ of the equation \eqref{eq:ham} 
is called a contractible orbit. The set  
$$ 
\Si_{H}:=\{\C A_H(x)\,|\, x  
\text{ is a contractible orbit }\} 
$$  
is called the action spectrum of $H$ and is compact 
\cite{HZ,S}. The following theorem of Schwarz is  fundamental~\cite{S}. 
 
\begin{theorem}\label{T:schwarz} 
Let $\Mo$ be a closed symplectically aspherical manifold. 
Let $H:S^1\times M\to \B R$ be a Hamiltonian whose time-$1$ map is not  
the identity. There are two contractible orbits $x,y\in M$ 
such that 
$$ 
-\int_0^1\max_MH_tdt\leq \C A_H(x)<\C A_H(y)\leq -\int_0^1\min_MH_tdt. 
$$ 
\end{theorem} 
 
Let $\psi \in \Ham\Mo$ be the time-1 map of a Hamiltonian 
$H:S^1\times M\to \B R$ and let $x,y\in M$ be its contractible orbits.  
Let $(\widetilde M,\widetilde \om) $ be the universal cover of $\Mo$ and 
$\widetilde \psi$ be a lift of $\psi$. Since $x,y$ are contractible 
orbits their lifts $\widetilde x,\widetilde y\in \widetilde M$ 
are contractible orbits as well (with respect to  
$\widetilde H:S^1\times \widetilde M\to \B R$ corresponding to $\widetilde  
\psi$). 
Take a curve $\ga:[0,1]\to \widetilde M$ with  
$\ga(0)=\widetilde x(0)$ and $\ga(1)=\widetilde y(0)$. 
Since $\widetilde M$ is simply connected there exists a disc  
$u:D^2\to \widetilde M$ with  boundary $\partial u = \psi\circ \ga - \ga$. 
Define 
$$ 
\Delta(\psi;x,y):= \int _{D^2}u^*(\om). 
$$ 
It is easy to see (using  symplectic asphericity and 
the contractibility of the orbits)  
that the definition does not depend on the 
choice of $\ga$ and $u$.  
 
The following lemma is proved by Polterovich~\cite{Po}. 
 
\begin{lemma} With the above notation the following hold: 
\begin{enumerate} 
\item 
$\Delta(\psi^n;x,y)=n\Delta(\psi;x,y)$; 
\item 
$\Delta (\psi;x,y)=\C A_H(y) - \C A_H(x).$ 
\end{enumerate} 
\end{lemma} 
 
As a simple corollary we get the following result. 
 
\begin{theorem} 
Let $\Mo$ be a symplectically aspherical manifold. 
Then the group $\Ham\Mo$ of Hamiltonian diffeomorphisms of $\Mo$ 
is torsion-free. 
\end{theorem}

\subsection{Hamiltonian representations of discrete groups} 
\label{SS:polterovich} 
 
Let $\Mo$ be a symplectically aspherical manifold. 
Let $\widetilde \om := p^*(\om)$, where $p:\widetilde M\to M$ 
is the universal cover. Clearly, $\widetilde \om$ is an exact 
two-form. Let $g$ be a Riemannian metric on $M$ and 
$\widetilde g$ the pull-back metric on the universal cover. Choose a point 
$x_0\in M$ and let $B(s)$ the ball of the radius $s$ centered at $x_0$. Let 
$u:\B R_+\to \B R_+$ be defined by 
\begin{equation}\label{eq:filling} 
u(s):=\inf_{d\al=\widetilde\om}\,\sup_{x\in B(s)}\,  
|\al(x)|_{\widetilde g}. 
\end{equation} 
The function $s\mapsto s\cdot u(s)$ is then strictly increasing. We 
define the {\em symplectic filling function} 
$v:\B R_+\to \B R_+$ to be its inverse. 
 
If $a_n$ and $b_n$ are positive sequences then we write 
$a_n\succeq b_n$ if there exists a constant $c>0$ such that 
$a_n\geq cb_n$ for all $n\in \B N$, and we write $a_n\sim b_n$ if 
$a_n\succeq b_n$ and $b_n\succeq a_n$. 
 
Let $G$ be a finitely generated group equipped with the word metric 
with respect to some finite generating set $A$. Let 
$\|g\|_A$ denotes the distance of $g\in G$ from the identity. 
The following result is proved in~\cite[Theorem 1.6.A]{Po}. 
 
\begin{theorem}\label{T:polt} 
Let $\Mo$ be a symplectically aspherical manifold. 
Let $A$ be a finite subset of $\Ham\Mo$ and $G$ be the subgroup of $\Ham\Mo$ 
generated by $A$. 
Then $\|g^n\|_A\succeq v(n)$ for all $g\in G$. 
\end{theorem} 
 
\begin{definition}\label{D:hyp} 
If the function $u: \R \to \R$ described in \eqref{eq:filling} is bounded, or 
equivalently $v\sim \id$,  
then the symplectic form $\omega$ is called {\em hyperbolic} and $\Mo$ is called 
{\em symplectically hyperbolic}. 
\end{definition} 
 
It is easy to see that boundedness of $u$ does not depend on the choice of the 
metric $g$ and the point $x_0$. 
 
K\"ahler manifolds of negative sectional curvature  
(e.g. closed surfaces of genus at least 2)  
are symplectically hyperbolic. The torus $T^2$ with the standard symplectic 
structure is not symplectically hyperbolic. 
 
The interest in symplectically hyperbolic manifolds is, in particular, motivated 
by the following fact. 
 
\begin{corollary}\label{C:polt} 
Let $\Mo$ be a symplectically hyperbolic manifold. Let $A$ be a finite subset of  
$\Ham\Mo$ and $G$ be the subgroup of $\Ham\Mo$ generated by $A$.  
Then every cyclic subgroup $\left<g\right>\subset G$ is undistorted with  
respect to the word metric given by $A$. 
In particular, $\Ham\Mo$ is torsion free. 
\end{corollary} 
 
\begin{remark} 
A cyclic subgroup $\left<g\right > \subset G$  
is undistorted if there 
exists a constant $C>0$ such that $|g^n|\geq C\cdot n$, 
for all $n\in \mathbb Z$. Here  
$|x|:=\min\{n\in \mathbb N\,|\, x=g_{i_1}^{p_1}\ldots g_{i_k}^{p_k}, \sum_i 
p_i=n\}$  
is the word norm of $x\in G$ with respect to a fixed finite set  
$\{g_1,\ldots g_m\}$ of generators of $G$. For example any 
cyclic subgroup of $\mathbb Z^n$ is undistorted. Finite cyclic subgroups 
are not undistorted. 
\end{remark} 
 
\subsection{More symplectically hyperbolic manifolds} 
\label{SS:more} 
 
Here we provide a source of examples of 
symplectically hyperbolic manifolds. The subject is 
treated in more detail in a forthcoming paper \cite{K}. 
Recall that a cohomology class is {\em bounded} if it is represented 
by a singular cochain whose values on singular 
simplices are uniformly bounded. Such cochains are 
called bounded as well.

\begin{lemma}\label{L:bounded} 
Let $\Mo$ be a symplectically aspherical manifold. 
If $\omega$ represents a bounded cohomology class 
then $\Mo$ is symplectically hyperbolic. 
\end{lemma} 
 
\begin{proof} 
We just sketch the main ideas of  proof referring to 
 \cite[Theorem 2.1]{K} for details. 
Let $p:\widetilde M \to M$ be the universal cover. 
We show that $\widetilde \om :=p^*(\om)=d\alpha$, 
where $\al$ is a form bounded with respect to 
the metric $\widetilde g$ induced from any metric $g$ 
on $M$. 
 
Let $c\in C^2(M;\R)$ be a bounded cochain representing 
the class $[\om]$. Since $\widetilde M$ is simply connected,  
we conclude that $p^*(c)= \delta (b)$ for some 
bounded real cochain $b\in C^1(\widetilde M)$.  
Let $K$ be a finite triangulation 
of $M$ and $K'$ the induced one of $\widetilde M$. Let  
$b',c'$ be simplicial cochains corresponding to $b$ and  
$c$ respectively. 
The standard construction of a differential form from a  
simplicial cochain~\cite[pages 148--149]{STh}) applied to $b'$ gives a 
bounded form $\al$ whose differential is $p^*(\om+d\beta)$ for 
some $\beta\in \Om^1(M)$. Hence $p^*(\om)=d(\al + p^*(\beta))$ while 
$\al+p^*(\beta)$ is clearly bounded. 
\end{proof}

\begin{corollary}\label{C:hyperbolic} 
Let $\Mo$ be a symplectically aspherical manifold.  
If $\pi_1(M)$ is hyperbolic then $\om$ 
is hyperbolic.  
In particular, if $M$ admits a Riemannian 
metric of negative sectional curvature then $\om $ is 
hyperbolic.  
\end{corollary} 

\begin{proof} 
By Proposition \ref{p:aspher}, 
$[\om]=f^*(\Om)$ for some $\Om\in H^2(\pi_1(M);\B R)$ and $f:M\to 
K(\pi_1(M),1)$.  
On the other hand, if 
$\pi_1(M)$ is hyperbolic then every cohomology class  
of $\pi_1(M)$ of degree greater than one is bounded.  
Hence $[\om]$ is bounded as well and we apply \lemref{L:bounded}. 
\end{proof} 
 Note that one can find the  definition of a hyperbolic group in \cite{Gr}.
\begin{corollary} 
Let $F$ be a closed oriented surface of genus at least $2$, 
and let $F\stackrel i\to M\stackrel p\to B$ be an oriented bundle over a  
surface $B$ of genus at least $1$. 
Then $M$ admits a hyperbolic symplectic form. 
\end{corollary} 
 
\begin{proof} 
Let $\om_B$ be an area form on $B$. Given any class $\Om\in H^2(M)$ 
with $i^*(\Om)\neq 0$, the Thurston construction~(\cite[Theorem 6.3]{MS1}) 
gives a symplectic form in the class $C\cdot p^*[\om_B]+ \Om$, where 
$C>0$ is a sufficiently large constant. 
 
Let $\Om$ be the Euler class of the 
bundle $V:=\ker dp \to M$ tangent to the fibers of $p$. 
According to Morita~\cite{Mo} this class is bounded, and therefore  
$p^*[\om_B]+\Om$ is. 
\end{proof} 
 
\subsection{The Ostrover trick \cite{Os}} 
\label{SS:ostrover} 
 
Let $B\subset M$ be an open subset and let $h\in \Ham\Mo$ be a 
Hamiltonian diffeomorphism such that $h(B)$ is disjoint from 
the closure of $B$. 
 
After slightly perturbing $h$  we may assume that its fixed 
points are all non-degenerate. 
Let $F:M\times [0,1]\to \B R$ be a normalized Hamiltonian  
such that $F(x,t)= C$ for some $C<0$ and $x\in M-B$, $t\in [0,1]$.  
Take $\psi_t=h\circ f_t$ where $f_t$ is the Hamiltonian  
flow generated by $F$. 
 
\begin{theorem}\label{T:ostrover} 
Let $\Mo$ be a symplectically aspherical manifold and 
let $\psi_t\in \Ham\Mo$ be as defined above. 
Then $\lim _{t\to \infty}d(\psi_t,\id)=\infty$. 
In particular, the Hofer diameter of $\Ham\Mo$ 
is infinite. 
\end{theorem} 
 
\begin{proof} 
Let $F_t$ be a Hamiltonian for $\psi_t$ and $H$ a one associated 
to $h$. 
According to Schwarz, the Hofer norm of a Hamiltonian 
diffeomorphism is bounded from below by the minimum 
of the action spectrum. That is we have that 
$$ 
d(\psi_t,\id) \geq \min\Si_{F_t}, 
$$ 
for all $t\in \B R$. 
Notice that $\psi_t$ and $h$ have the same fixed points. 
Moreover, Ostrover proved (Proposition 2.6. in \cite{Os}) that 
$$ 
\C A_{F_t}(x)=\C A_H(x) - t\cdot C. 
$$ 
Combining this with the previous inequality we get 
that the Hofer norm of $\psi_t$ tends to infinity as $t$ does. 
\end{proof} 
 
\section{Application of symplectic asphericity to circle actions on symplectic  
manifolds} 
 
If  a closed symplectic manifold $(M,\omega)$ admits a Lie group action  
preserving $\omega$, the manifold must satisfy various topological  
restrictions. The nature of these restrictions is to some extent understood. It  
comes from  Morse theory. For example, assume that $(M,\omega)$ admits a  
circle action that is {\it Hamiltonian}, i.e. $i_X\omega=d\mu$ for some smooth  
function $\mu: M\rightarrow \mathbb{R}$ and the fundamental vector field $X$  
determined by the circle action. This function $\mu$ is the {\it  
Bott-Morse} function, and this forces certain restrictions on the  
topology of $M$. Finding these conditions is now a huge research area.  
Restrictions on the equivariant cohomology of $M$ are given in \cite{Ki, TW},  
 on characteristic classes, signature and Novikov numbers can be  
found in \cite{Fa, Fe}, restrictions on the topology of orbits in \cite{Oz,Ko},  
on Massey products in \cite{ST}. One can work in   pure homotopic setting  
of {\it cohomologically symplectic manifolds} with circle actions, and  
still get non-trivial restrictions on the equivariant cohomology of $M$ and the  
set of fixed points \cite{A1,A2}. Many results in this theory are obtained as  
variations of the following fundamental fact. If $G$ is a torus acting in a  
Hamiltonian way on a closed  symplectic manifold $M$, then the fiber bundle 
$$ 
\CD 
M @>{i}>> EG\times _GM @>>> BG 
\endCD 
$$ 
is {\it totally non-cohomologous to zero}, that is,  
$$ 
i^*: H^*_G(M;\R):=H^*(EG\times _GM;\R)\rightarrow H^*(M;\R) 
$$ 
is onto~\cite{Ki}. In the sequel we will call this the TNCZ property. In 
general, TNCZ does not hold for cohomologically symplectic manifolds with circle  
actions \cite{A1}. However, some properties of symplectic manifolds with  circle 
actions   
do have cohomologically symplectic analogues. Such results  follow  
as combinations of the localization theorem for equivariant cohomology with  
various additional algebraic assumptions. For example, if one {\it imposes}  
TNCZ condition on a circle action on a   
cohomologically symplectic manifold, one obtains that the set of fixed points of 
this action has at least   
two components, as in the case of true symplectic actions \cite{A1,A2}(there are 
some additional technical assumptions, but we   
don't discuss them). 
 
However, it seems that the natural boundaries of the theory are still not  
explored. For example, we 
don't know any examples of closed symplectic manifolds endowed with circle  
actions but whose topological properties differ from the ones established for  
symplectic manifolds with {\it symplectic} circle actions.  
 
If one imposes the condition of symplectic asphericity, more restrictions can  
be found. 
In particular, Ono~\cite{O} found  restrictions on the fundamental group of $M$  
in the presence of the symplectic asphericity condition. On the other hand,  
results of Lupton and Oprea \cite{LO} suggest that these restrictions may have a 
purely homotopic  
nature.  
 
\subsection{Symplectic asphericity as an obstruction} 
Below we shall show how  symplectic asphericity obstructs the existence of 
symplectic circle  
actions. 
 
 
The following result is a weak version of~~\cite[Theorem 4.16]{LO} (see also  
~~\cite{Op}). 
  
\begin{theorem}\label{t:lo} Any $S^1$-action on a symplectically aspherical  
manifold has no fixed points and hence is not Hamiltonian. 
\end{theorem}

\begin{theorem}\label{t:ono}  
Let $(M,\omega)$ be a closed symplectically aspherical manifold. If each  
Abelian subgroup of $\pi_1(M)$ is cyclic, then $M$ does not admit  non-trivial  
symplectic  circle actions. 
\end{theorem} 
 
\begin{proof} First, it follows from the symplectic asphericity assumption and 
the previous theorem, that a circle  
action on $M$ cannot be Hamiltonian. However,  
the action also cannot be non-Hamiltonian, which follows from the following  
argument. By way of contradiction assume that there exists a non-Hamiltonian  
action $a: S^1\to \Symp(\Mo)$ and let $X$ denote the vector field of this  
action. Since $a$ is non-Hamiltonian, we conclude that $[i_X \omega]\ne 0$ in  
$H^1(M;\R)$. Hence there exists a loop $A: S^1\to M$ such that $\la  
[i_X\omega,[A]\ra\ne 0$ where $[A]$ is the homology class of the loop $A$. Now  
consider the map 
\[ 
a(A): T^2\to M,\quad a(A)(s,t)=a(s)(A(t)), 
\] 
and it is easy to see that 
\[ 
\la a(A)^*[\omega], [T^2]\ra=\int_{T^2}a(A)^*\omega \ne 0. 
\] 
 
Since the image of the homomorphism $a(A)_*:\pi_1(T^2)\to \pi_1(M)$ is cyclic,  
we conclude that there exists a simple loop $C$ on $T^2$ whose homotopy class  
is nontrivial and belongs to the kernel of $a(A)_*$.  
Hence  
we 
get a map $f: S^2 \to M$ with $\la  
[\omega], f_*[S^2]\ra \ne 0$. But this contradicts the symplectic  
asphericity of $\Mo$. 
\end{proof} 
   
\begin{remark} \hfill 
\begin{enumerate}  
\item   
Theorem 8.3 was proved by Ono in \cite{O}. 
\item  
 It would be nice to find other  
topological restrictions on symplectically aspherical manifolds with circle  
actions, which follow from the condition that the action is symplectic. It is 
conceivable, that Ono's condition is stronger: probably, it implies the 
non-existence of {\it smooth} circle actions on symplectically aspherical 
manifolds. A possible proof would go along the following lines. There is a 
notion of cohomologically Hamiltonian circle action \cite{LO}. One could try to 
use it instead of the condition $[i_X\omega]=0$. However, we did not work out 
the details.  
\item 
In fact, \cite[Theorem 4.16]{LO} is a purely cohomological (and hence more 
general) version of \theoref{t:lo}. 
\end{enumerate} 
\end{remark} 
 
\subsection{On a problem of Taubes and related questions} The following problem  
was posed by Taubes~\cite{Ba},~\cite{FV}. 
 
\begin{question} Assume that a $4$-manifold of the form $M^4=N^3\times S^1$ is  
symplectic. Is it true 
that $N$ fibers over $S^1$? 
\end{question} 
 
Note that if $N$ fibers over $S^1$ then $M$ admits a symplectic  
structure~\cite{Ba}. The question was answered in the affirmative in several  
important cases  (see~\cite{FV}), but in general is still open. One can  
reformulate it in the following form. 
 
\begin{question} Assume that $(M^4,\omega)$ is a closed symplectic manifold  
endowed with a free circle action. Does it admit a circle action preserving the  
given symplectic form $\omega$? 
\end{question} 
 
It seems interesting to ask the similar question for manifolds of arbitrary  
dimension (which might be easier, since more classical differential topology  
methods are available). Note that the symplectic asphericity condition might 
help  
in looking for counterexamples in higher dimensions. Ono's theorem yields a  
justification: we know that any circle action on a symplectically aspherical  
$(M,\omega)$  cannot have fixed points, and we have a restriction on the  
fundamental group. 
\vskip6pt  
{\bf Acknowledgment.} The research of the second author was supported by NSF 
grant 0406311. The third author was supported by the Ministry of Science and 
Higher Education of Poland, research project no. 1P03A 03330.

\end{document}